\documentclass[12pt]{article}
\usepackage{amssymb,amsmath,amsthm}
\usepackage{makeidx,graphicx,amssymb,amsmath}
\usepackage{amsthm,fullpage, wrapfig}
\headsep=0mm \numberwithin{equation}{section}

\title{\bf \large The Rayleigh-Taylor instability  for the Musket
 free boundary problem.}

\author{O. V. Galtsev, A. M. Meirmanov}
\date{}

\theoremstyle{plain}

\theoremstyle{definition}

\theoremstyle{remark}

\numberwithin{equation}{section}
\renewcommand{\mathbf}[1]{\mbox{\boldmath$#1$}}

\begin{document}

\maketitle \small
\begin{abstract}
 The present paper is devoted to the joint
motion of two immiscible incompressible liquids  in porous media.
The liquids have different densities and initially separated by a
surface of strong discontinuity (free boundary). We discuss the results
of numerical simulations for exact free boundary problems on the
microscopic level for the absolutely rigid solid skeleton and for
the elastic solid skeleton of different geometries. The problems
have a natural small parameter, which is the ratio of average pore
size to the size of the domain in consideration. The formal limits
as $\varepsilon\searrow 0$  results homogenized models, which are
the Muskat problem in the case of the absolutely rigid solid
skeleton, and the viscoelastic Muskat problem in the case of the
elastic solid skeleton.  The last model preserves a free boundary
during the motion, while in the first model instead of the free
boundary appears a mushy region, occupied by
a mixture of two fluids.\\
\end{abstract}

\noindent \textbf{Key words:} Musket
  problem, the Rayleigh-Taylor instability,
Stokes and Lam\'{e}'s equations.\\

\noindent \textbf{MOS subject classification:} 35R35, 35M13, 35B27\\

\normalsize

\addtocounter{section}{1} \setcounter{equation}{0}
\begin{center} \textbf{ \S 1. Introduction}
\end{center}
The Rayleigh–Taylor instability, or RT instability (after Lord Rayleigh \cite{LR} and G. I. Taylor \cite{GT}), is an instability of an interface between two fluids of different densities, which occurs when the lighter fluid is pushing the heavier fluid.  The equivalent situation occurs when gravity is acting on two fluids of different density — with the dense fluid above a fluid of lesser density. As the instability develops, downward-moving irregularities are quickly magnified into sets of inter-penetrating Rayleigh–Taylor fingers. Therefore RT instability is sometimes qualified to be a fingering instability \cite{CHP}. In the present paper we consider RT instability in an inhomogeneous fluid filling the pores in the solid skeleton. In addition to the undoubted theoretical interest, this problem  is important for a number of  practical problems.
For example, the description of the displacement of one liquid by another in porous
media. This problem still needs correct mathematical model.

There are different types of mathematical models, but we are interested only in some of the fundamental models of continuum mechanics (such as, for example, Stokes equations for a slow motion of a viscous liquid, or Lam\'{e}'s equations for displacements of an elastic solid body), or in models asymptotically close to above mentioned ones.

Among mathematical models of a joint motion of two immiscible liquids the most trustable (or physically correct) one is the Muskat problem, suggested by M. Muskat\cite{MM}. This model describes a filtration of two immiscible incompressible liquids of different viscosities and different densities, divided by some moving boundary (free boundary). The  motion of the first  liquid under the gravity in the domain $\Omega^{+}(t)$ with a constant viscosity $\mu$ and a constant density $\rho_{f}^{+}$ is governed by the Darcy system of filtration
\begin{equation}
\boldsymbol{v}^{+}=-\frac{k}{\mu}\nabla p_{f}^{+}+\rho_{f}^{+}\boldsymbol{e}_{2},\quad
\nabla\cdot\,\boldsymbol{v}^{+}=0,
\quad \boldsymbol{x}\in \Omega^{+}(t),\label{1.1}
\end{equation}
for the macroscopic velocity $\boldsymbol{v}^{+}$ and the pressure $p_{f}^{+}$ of the liquid. The unit vector $\boldsymbol{e}_{2}$ coincides with the direction of the gravity.

Correspondingly, the motion of the second  liquid under the gravity in the domain $\Omega^{-}(t)$ with a constant viscosity $\mu$ and a constant density $\rho_{f}^{-}$  is governed by  the Darcy system of filtration
\begin{equation}
\boldsymbol{v}^{-}=-\frac{k}{\mu}\nabla p_{f}^{-}+\rho_{f}^{-}\boldsymbol{e}_{2},\quad
\nabla\cdot\,\boldsymbol{v}^{-}=0,
\quad \boldsymbol{x}\in \Omega^{-}(t),\label{1.2}
\end{equation}
for the macroscopic velocity $\boldsymbol{v}^{-}$ and the pressure $p_{f}^{-}$. On the common free boundary $\Gamma(t)=\partial\Omega^{+}(t)\bigcap\partial\Omega^{-}(t)$ pressures and normal velocities are continuous:
\begin{equation}
p_{f}^{+}=p_{f}^{-},\quad \boldsymbol{x}\in \Gamma(t),\label{1.3}
\end{equation}
\begin{equation}
\boldsymbol{v}^{+}\cdot \boldsymbol{n}=
\boldsymbol{v}^{-}\cdot \boldsymbol{n}=V_{n},
\quad \boldsymbol{x}\in \Gamma(t),\label{1.4}
\end{equation}

\begin{wrapfigure}{l}{90pt}
\includegraphics[width=35mm,height=40mm]{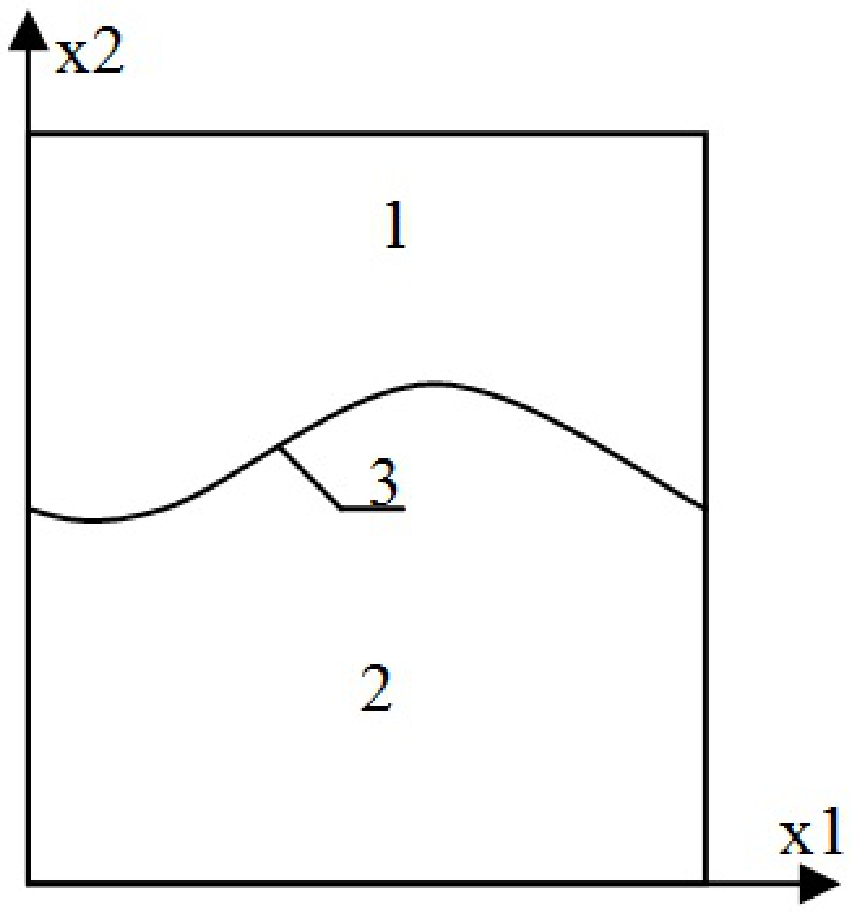}
\phantom{.}\vskip 0 cm $\phantom{.}$ {\scriptsize  Fig. 1.
1 -- domain $\Omega^{+}(t)$, 2 -- domain
$\Omega^{-}(t)$, 3 -- free boundary  $\Gamma(t)$.}
\end{wrapfigure}

where $\boldsymbol{n}$ is the unit normal vector to the boundary
$\Gamma(t)$ at the point $\boldsymbol{x}\in \Gamma(t)$ and $V_{n}$
is the velocity in the normal direction  of the boundary $\Gamma(t)$ at
the point $\boldsymbol{x}\in \Gamma(t)$.

Condition (\ref{1.4}) means that the  boundary $\Gamma(t)$
is a material surface. That is, it consists of the same set of material points during the motion. This fact permits  a weak formulation of
the Muskat problem. So, we define the pressure $p_{f}$ of the
inhomogeneous liquid as
\[
p_{f}=p_{f}^{+}\,\, \mbox{if}\,\,\boldsymbol{x}\in \Omega^{+}(t),
\,\,\,p_{f}=p_{f}^{-}\,\, \mbox{if}\,\,
\boldsymbol{x}\in \Omega^{-}(t),
\]
the density $\rho_{f}$ as
\[
\rho_{f}=\rho_{f}^{+}\,\, \mbox{if}\,\,\boldsymbol{x}\in \Omega^{+}(t),
\,\,\,\rho_{f}=\rho_{f}^{-}\,\, \mbox{if}\,\,
\boldsymbol{x}\in \Omega^{-}(t),
\]
and the velocity $\boldsymbol{v}$ as
\[
\boldsymbol{v}=\boldsymbol{v}^{+}\,\, \mbox{if}\,\,\boldsymbol{x}\in \Omega^{+}(t),
\,\,\,\boldsymbol{v}=\boldsymbol{v}^{-}\,\, \mbox{if}\,\,
\boldsymbol{x}\in \Omega^{-}(t).
\]
Then the unknown functions $\boldsymbol{v}$, $p_{f}$, and $\rho_{f}$ satisfy the Darcy system of filtration in the form
\begin{equation}
\boldsymbol{v}=-\frac{k}{\mu}\nabla p_{f}+\rho_{f}\,\boldsymbol{e}_{2},\quad
\nabla\cdot\,\boldsymbol{v}=0,
\quad\boldsymbol{x}\in \Omega,\quad t>0,\label{1.5}
\end{equation}
and transport equations
\begin{equation}
\frac{d\rho_{f}}{dt}\equiv\frac{\partial\rho_{f}}{\partial t}+
\boldsymbol{v}\cdot\nabla\,\rho_{f}=0,
\quad\boldsymbol{x}\in \Omega,\quad t>0,\label{1.6}
\end{equation}
The first equation in (\ref{1.5}) (Darcy law) is understood in an usual sense almost everywhere in $\Omega_{T}=\Omega\times (0,T)$, and the second equation (continuity equation) is understood  in the sense of distribution. Transport equation is understood in a sense of distributions, if we use  equalities
\[
\boldsymbol{v}\cdot\nabla\mu=\nabla\cdot(\boldsymbol{v}\mu),\quad
\boldsymbol{v}\cdot\nabla\,\rho_{f}=\nabla\cdot(\boldsymbol{v}\rho_{f}).
\]

The problem is endowed with homogeneous boundary condition
\begin{equation}
\boldsymbol{v}\cdot\boldsymbol{n}=0,
\quad\boldsymbol{x}\in S=\partial\Omega,\quad t>0,\label{1.7}
\end{equation}
where $\boldsymbol{n}$ is the normal vector to the boundary $S$, and
initial conditions
\begin{equation}
\rho_{f}(\boldsymbol{x},0)=\rho_{f}^{0}(\boldsymbol{x}),
\quad\boldsymbol{x}\in \Omega,\label{1.8}
\end{equation}
with discontinuous initial data:
\[
\rho_{f}^{0}(\boldsymbol{x})=\rho_{f}^{+},\,\,\mbox{if}\,\,
\boldsymbol{x}\in \Omega^{+}\,\,\mbox{and}\,\,
\rho_{f}^{0}(\boldsymbol{x})=\rho_{f}^{-},\,\,\mbox{if}\,\,
\boldsymbol{x}\in \Omega^{-}.
\]
So, one has two settings of the same Muskat problem.
In both cases the problem is easy  to formulate, but almost impossible to solve. For this reason, very little is known either about  classical solutions or on weak solutions. There are only few results on a classical solvability locally in time or globally in time, but near explicit solutions, and there is no any result on a weak solvability (see,  Ref.~\cite{FHY}, Ref.~\cite{ER},  Ref.~\cite{SCH} and references there).

Following R. Burridge and J. Keller Ref.~\cite{BK}, let us try to find more general physically correct mathematical models describing the same physical process. To explain idea we, at first,consider the Darcy system of filtration, which is responsible for the dynamics in the Muskat problem. It is well-known that this system is an asymptotic limit of the Stokes system for an incompressible viscous liquid, when dimensionless pore size  tends to zero (see \cite{BK}, \cite{SP}). R. Burridge and J. Keller have suggested to consider more general system
\begin{equation}
\alpha_{\tau}\rho^{\varepsilon}
\frac{\partial^2\boldsymbol w}{\partial t^2}=
\nabla\cdot\big(\chi^{\varepsilon}\alpha_{\mu}
\mathbb{D}(x,\frac{\partial\boldsymbol w}{\partial t})+(1-\chi^{\varepsilon})\alpha_{\lambda}
\mathbb{D}(x,\boldsymbol{w})-p\mathbb{I}\big)+
\rho^{\varepsilon}\,\boldsymbol{e}_{2},
\label{1.9}
\end{equation}
\begin{equation}
p+\alpha_{p}\nabla\cdot\,\boldsymbol{w}=0,\label{1.10}
\end{equation}
for the displacement $\boldsymbol{w}$ and pressure $p$ of the continuum medium. The microscopic system (\ref{1.9}), (\ref{1.10}) describes the joint motion of the viscous liquid in pore space  and elastic solid skeleton and is understood in the sense of distributions. Roughly speaking, this system contains the Stokes system for the viscous liquid in the pore space $\Omega_{f}$, the Lam\'{e}'s system for the solid skeleton in $\Omega_{s}$ and boundary condition (the continuity of the normal stresses) on the common boundary $\Gamma=\partial\Omega_{f}\cap\partial\Omega_{s}$.  In (\ref{1.9}) $\mathbb{D}(x,\boldsymbol{w})$ is the symmetric part of $\nabla\boldsymbol{w}$,  $\chi^{\varepsilon}$ is the characteristic function of the pore space $\Omega_{f}$, $\varepsilon=l/L$ is the dimensionless pore size,
\[
\alpha_{\tau}=\frac{ L}{g \tau^2}, \quad \alpha_{\mu}=\frac{2\mu
}{\tau Lg\rho_{0}},\quad \alpha_{\lambda}=\frac{2 \lambda}{Lg\rho_{0}},
\]
\[
\rho^{\varepsilon}=\rho_{f}\chi^{\varepsilon}+\rho_{s}(1-\chi^{\varepsilon}),\quad
\alpha_{p}=\alpha_{p,f}^{2}\chi^{\varepsilon}+
\alpha_{p,s}^{2}(1-\chi^{\varepsilon}),
\]
$l$ is an average size of pore,  $L$ is a characteristic size of the domain in consideration, $\tau$ is a characteristic time of the process, $\rho_{f}$ and $\rho_{s}$ are the respective mean dimensionless densities of the liquid in pores and the solid skeleton correlated with the mean density of water $\rho_0$,  $\alpha_{p,f}$ and $\alpha_{p,s}$ are the respective dimensionless speed of sound of the liquid in pores and the solid skeleton,  $g$ is the value of acceleration of gravity,  $\mu$  is the viscosity of fluid, and $\lambda$ is the elastic Lam\'{e}'s constant. In what follows we suppose the structure of the domains $\Omega_{f}$ and $\Omega_{s}$ is periodic with a period $\varepsilon$. That is $\chi^{\varepsilon}(\boldsymbol{x})=\chi(\boldsymbol{x}/\varepsilon)$ with 1-- periodic function $\chi(\boldsymbol{y})$.

Theoretically the system (\ref{1.9}), (\ref{1.10}) with corresponding boundary and initial conditions most accurately describes the given physical process, but it still has no any practical significance. This one appears only after homogenization. So, we have to let all dimensionless criteria $\alpha_{\tau},\,\alpha_{\mu}$ and $\alpha_{\lambda}$ to be variable functions, depending on the small parameter $\varepsilon$, and find all limiting regimes as $\varepsilon\rightarrow 0$.

It is clear that limiting regimes of the system (\ref{1.9}), (\ref{1.10}) depend on the dimensionless criteria
\[
\tau_{0}=\lim_{\varepsilon\searrow 0} \alpha_{\tau},
\quad \mu_{0}=\lim_{\varepsilon\searrow 0} \alpha_{\mu},\quad
\lambda_{0}=\lim_{\varepsilon\searrow 0} \alpha_{\lambda},
\]
\[
c_{f}=\lim_{\varepsilon\searrow 0}\alpha_{p,f}
\quad c_{s}=\lim_{\varepsilon\searrow 0}\alpha_{p,s},\quad
\mu_{1}=\lim_{\varepsilon\searrow 0}\frac{\alpha_{\mu}}{\varepsilon^{2}}
\]
(see \cite{AM}). For example, for filtration $\tau_{0}=0$, for absolutely rigid solid skeleton $\lambda_{0}=\infty$, and for incompressible media $c_{f}=c_{s}=\infty$. Therefore, different sets of criteria lead to different homogenized models, asymptotically closed to the basic one. All of these models describe the same physical process, but with a different degree of approximation.

For a filtration of an incompressible liquid ( $c_{f}=\infty$) in an incompressible solid skeleton ( $c_{s}=\infty$)  the standard homogenized model for fixed  $\alpha_{\mu}=\mu_{0}$, $\alpha_{\lambda}=\lambda_{0}$, and $\alpha_{\tau}=\tau_{0}$ has a form
\[
\tau_{0}\hat{\varrho}\frac{\partial^{2}\boldsymbol{w}}{\partial t^{2}}=\nabla\cdot\,\widehat{\mathbb{P}}+\hat{\varrho}\,\boldsymbol{e}_{2},\,\,
\hat{\varrho}=m\rho_{f}+(1-m)\rho_{s},
\]
\[
\widehat{\mathbb{P}}=
\mu_{0}\mathfrak{N}_{1}:\mathbb{D}(x,\frac{\partial\boldsymbol{w}}{\partial
t})+ \lambda_{0}\mathfrak{N}_{2}:\mathbb{D}(x,\boldsymbol{w})+
\int_{0}^{t}\mathfrak{N}_{3}(t-\tau):\mathbb{D}(x,\boldsymbol{w}(\boldsymbol{x},\tau))d\tau.
\]
But the same system (\ref{1.9}), (\ref{1.10}) on the microscopic level has another asymptotic limits, like Darcy system of filtration ($\tau_{0}=0$, $c_{f}=c_{s}=\infty$, $\lambda_{0}=\infty$, $\mu_{0}=0$, $0<\mu_{1}<\infty$):
\[
\frac{\partial\boldsymbol{w}}{\partial t}=\frac{1}{\mu_{1}}
\mathbb{B}^{f}\cdot(-\nabla p+\rho_{f}\,\boldsymbol{e}_{2}),\,\,
\nabla\cdot\,\frac{\partial\boldsymbol{w}}{\partial t}=0,
\]
or  Terzaghi -- Biot system of poroelasticity ($\tau_{0}=0$, $c_{f}=c_{s}=\infty$, $0<\lambda_{0}<\infty$, $\mu_{0}=0$, $0<\mu_{1}<\infty$):
\[
\frac{\partial\boldsymbol{w}}{\partial t}=
m\frac{\partial\boldsymbol{u}}{\partial t}+
\frac{1}{\mu_{1}} \mathbb{B}^{f}\cdot(-\nabla p+\rho_{f}\,\boldsymbol{e}_{2}),
\]
\[
\nabla\cdot\,\big(\frac{\partial\boldsymbol{w}}{\partial t}+
(1-m)\frac{\partial\boldsymbol{u}}{\partial t}\big)=0,
\]
\[
\nabla\cdot\,\big(\lambda_{0}\mathfrak{N}_{0}:\mathbb{D}(x,\boldsymbol{u})\big)-
\nabla p+\hat{\varrho}\,\boldsymbol{e}_{2}=0,
\]
or the system of viscoelastic filtration ($\tau_{0}=0$, $c_{f}=c_{s}=\infty$, $0<\lambda_{0},\,\mu_{0}<\infty$):
\[
\nabla\cdot\,\widetilde{\mathbb{P}}-\nabla p+\hat{\varrho}\,\boldsymbol{e}_{2}=0,\,\,
\nabla\cdot\,\boldsymbol{w}=0,
\]
\[
\widetilde{\mathbb{P}}=\mu_{0}\mathfrak{N}_{4}:
\mathbb{D}(x,\frac{\partial\boldsymbol{w}}{\partial t})+
\lambda_{0}\mathfrak{N}_{5}:\mathbb{D}(x,\boldsymbol{w})+
\int_{0}^{t}\mathfrak{N}_{6}(t-\tau):\mathbb{D}(x,\boldsymbol{w}(\boldsymbol{x},\tau))d\tau
\]
(see  \cite{AM},  \cite{KT}, \cite{MB}).

The same scheme we may apply to the free boundary problem. On the microscopic level this model for filtration in incompressible media has the form
\begin{equation}
\nabla\cdot\big(\chi^{\varepsilon}\alpha_{\mu}
\mathbb{D}(x,\frac{\partial\boldsymbol {w}^{\varepsilon}}{\partial t})+
(1-\chi^{\varepsilon})\alpha_{\lambda}
\mathbb{D}(x,\boldsymbol{w}^{\varepsilon})-p^{\varepsilon}\mathbb{I}\big)+
\rho^{\varepsilon}\,\boldsymbol{e}_{2}=0,
\label{1.11}
\end{equation}
\begin{equation}
\nabla\cdot\,\boldsymbol{w}^{\varepsilon}=0,\label{1.12}
\end{equation}
\begin{equation}
\frac{d\rho^{\varepsilon}}{dt}\equiv\frac{\partial\rho^{\varepsilon}}{\partial t}+
\frac{\partial\boldsymbol {w}^{\varepsilon}}{\partial t}
\nabla\rho^{\varepsilon}=0.\label{1.13}
\end{equation}
This problem has been completely studied in \cite{AM1}. The main result there states, that for any $\varepsilon >0$ there exists a weak solution $\{\boldsymbol{w}^{\varepsilon},\,p^{\varepsilon},\,\rho^{\varepsilon}\}$ to the free boundary problem  (\ref{1.11}) -- (\ref{1.13}) and for $\alpha_{\mu}=\mu_{0}$, $\alpha_{\lambda}=\lambda_{0}$ the corresponding sequences converge as $\varepsilon \searrow 0$ to the solution $\{\boldsymbol{w},\,p,\,\rho\}$ of the (homogenized) Muskat problem for the viscoelastic filtration
\begin{equation}
\nabla\cdot\,\widetilde{\mathbb{P}}-\nabla p+\rho\,\boldsymbol{e}_{2}=0,\,\,
\nabla\cdot\,\boldsymbol{w}=0,\,\,\frac{d\rho}{dt}=0.\label{1.14}
\end{equation}
The proof of this result has essentially used the notion of the two -- scale convergence \cite{NGU}.

On the other hand the formal limit in (\ref{1.11}) -- (\ref{1.13}) as
$\varepsilon \searrow 0$ for $\alpha_{\lambda}\rightarrow\infty$  and  $\alpha_{\mu}=\mu_{1}\varepsilon^{2}$ leads to the Muskat problem (\ref{1.5}), (\ref{1.6}). In is clear, that the same formal limit we obtain if as a basic model on the microscopic level we consider the free boundary problem for the Stokes system
\begin{equation}
\nabla\cdot\big(\mu_{1}\varepsilon^{2}
\mathbb{D}(x,\frac{\partial\boldsymbol{w}^{\varepsilon}}{\partial t})-p^{\varepsilon}\mathbb{I}\big)+\rho^{\varepsilon}\,\boldsymbol{e}_{2}=0,\,\,
\boldsymbol{x}\in \Omega^{\varepsilon}_{f},\,t>0,
\label{1.15}
\end{equation}
\begin{equation}
\nabla\cdot\,\boldsymbol{w}^{\varepsilon}=0,\,\,
\frac{d\rho^{\varepsilon}}{dt}=0,\,\,\boldsymbol{x}\in \Omega^{\varepsilon}_{f},\,t>0,
\,\,\boldsymbol{w}^{\varepsilon}=0,\,
\boldsymbol{x}\in \partial\Omega^{\varepsilon}_{f},\,t>0,
\label{1.16}
\end{equation}
only in the pore space $\Omega^{\varepsilon}_{f}$. The last problem has been  studied in \cite{AMY}, where authors proved that for any $\varepsilon >0$ the problem (\ref{1.15}) has the unique classical solution $\{\boldsymbol{w}^{\varepsilon},\,p^{\varepsilon},\,\rho^{\varepsilon}\}$.

The goal of the present paper is numerical simulations for the problem (\ref{1.11}) -- (\ref{1.13}) with $\alpha_{\mu}=\mu_{0}$, $\alpha_{\lambda}=\lambda_{0}$, and for the problem (\ref{1.15}), (\ref{1.16}). We have done numerical results for two different structures of the pore space: disconnected capillaries and disconnected solid skeleton in the unit square in $\mathbb{R}^{2}$.

Numerical simulations of the problem in a single capillary in the
absolutely rigid skeleton show the coincidence with  results of
\cite{BJG}. On the Figure 2 we may see the smooth free boundary (the
surface of a strong discontinuity) in the capillary at different
times.\\
\\
\\
\\

\begin{center}
\includegraphics[width=60mm,height=35mm]{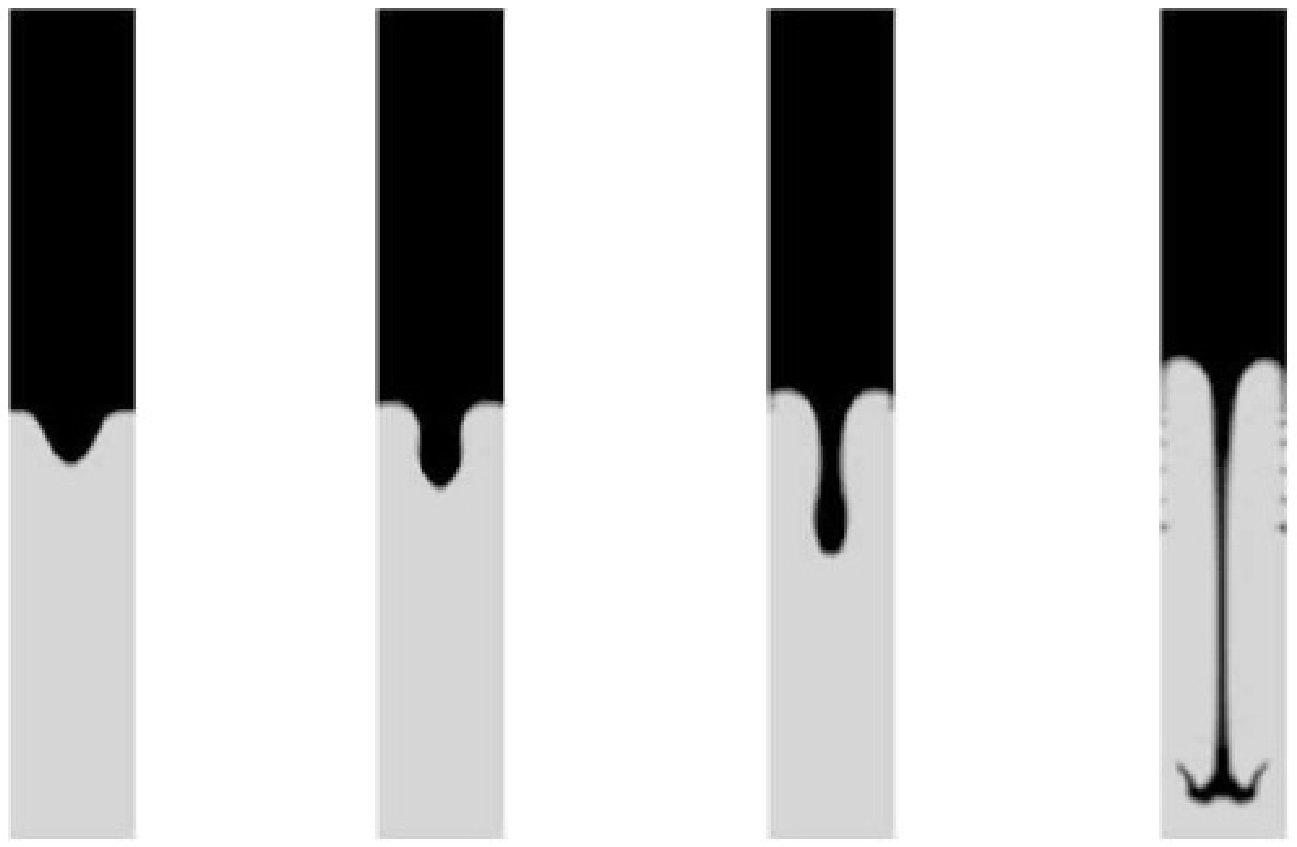}
\end{center}
\begin{center}
\ \phantom{.}\hskip - 1mm $\scriptstyle t = 0.003$\hskip 6mm
$\scriptstyle t = 0.255$\hskip 6mm $\scriptstyle t = 0.525$\hskip
6mm  $\scriptstyle t = 1.035$
\end{center}
\begin{center}
{\scriptsize Fig. 2. The surface of a strong discontinuity at
different times.}
\end{center}

The numerical results show that the motion of the liquids describing by (\ref{1.15}), (\ref{1.16}) is affected by three factors: the ratio $\delta = \rho^{+}/\rho^{-}$ of the densities $\rho^{+}$ and $\rho^{-}$ of the liquids on the top and on the bottom, the viscosity $\mu$ of fluids and the pore size $\varepsilon$. Changing these parameters one gets different scenarios of Rayleigh-Taylor instability. For the problem (\ref{1.11}) -- (\ref{1.13}) the process is affected by the same parameters $\delta$, $\mu$, $\varepsilon$ as before, and by the additional parameter $\lambda$, which is the elastic Lam\'{e}'s constant.

The limiting procedure ($\varepsilon \searrow 0$) is modeled by increasing the number of capillaries for the first geometry, and the number of elementary squares for the second geometry. Therefore, we may assume that for sufficiently small $\varepsilon$ the problem (\ref{1.15}), (\ref{1.16}) describes the classical Muskat problem, while the problem (\ref{1.11}) -- (\ref{1.13}) describes the Muskat problem for viscoelastic filtration.

Figure 3 shows numerical simulations for the first geometry for the model (\ref{1.15}), (\ref{1.16}) of the motion in an absolutely rigid solid skeleton (above), and for the model (\ref{1.11}) -- (\ref{1.13}) of the motion in the elastic solid skeleton (below).

\begin{center}
\includegraphics[width=165mm,height=60mm]{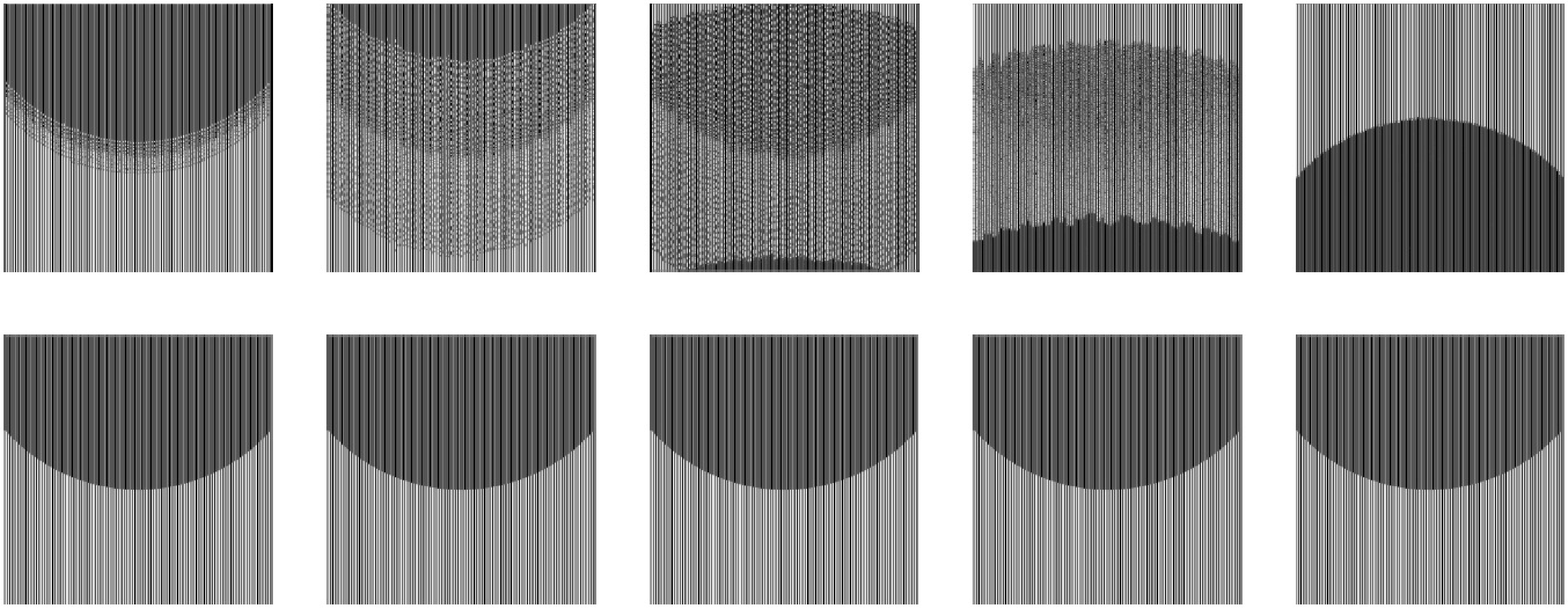}
\vskip - 0.1 cm \phantom{.} \hskip 0mm $\scriptstyle t = 50$\hskip
28mm $\scriptstyle t = 860$ \hskip 22mm $\scriptstyle t = 2631$
\hskip 24mm $\scriptstyle t = 3012$ \hskip 22mm $\scriptstyle t =
4873$
\end{center}
\begin{center}
{\scriptsize Fig. 3.  Disconnected capillaries: numerical simulation for the absolutely rigid (above) and for the elastic solid skeleton.}
\end{center}

\begin{center}
\includegraphics[width=70mm,height=70mm]{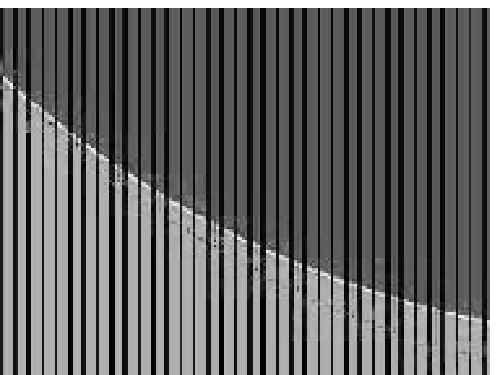}
\includegraphics[width=70mm,height=70mm]{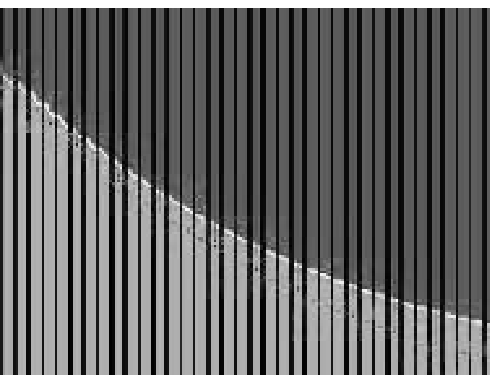}

\vskip 0 cm \phantom{.} \hskip 0mm $\scriptstyle t = 50$\hskip 65mm
$\scriptstyle t = 860$
\end{center}
\begin{center}
\includegraphics[width=70mm,height=70mm]{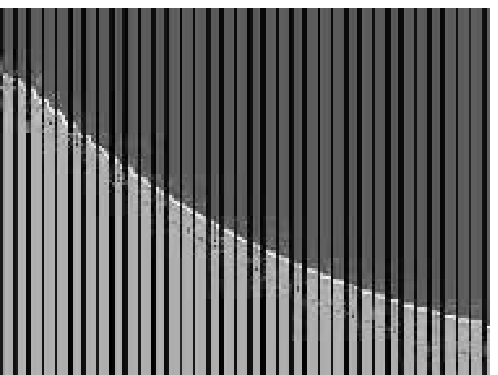}
\includegraphics[width=70mm,height=70mm]{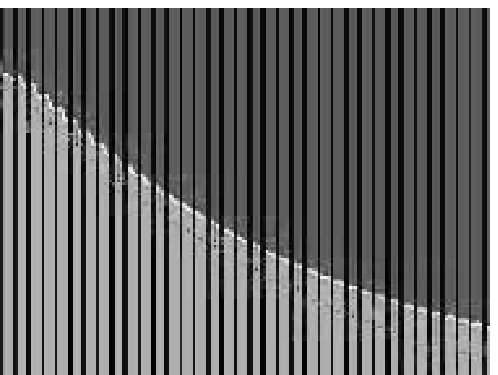}
\vskip 0 cm \phantom{.} \hskip 0mm $\scriptstyle t = 2631$ \hskip
60mm $\scriptstyle t = 3012$
\end{center}

\begin{center}
{\scriptsize Fig. 4. The case of an elastic solid skeleton in a large
scale.}
\end{center}
We also  carried out numerical simulations for different values of  $\lambda$
and $\delta$:

\begin{itemize}
 \item for $\varepsilon = 2*10^{-5}$, $\lambda \searrow 0$ and $\delta = 1.25$, there is a changing of the interface of fluids;
 \item for $\varepsilon = 2*10^{-5}$, $\lambda = 0.5$ and $\delta \searrow \infty$, there is a changing of the interface of fluids;
 \item for $\varepsilon = 2*10^{-5}$, $\lambda = 0.5 $ and $\delta \searrow 1$, there is no a changing of the interface of fluids.
\end{itemize}

The same conclusion is valid for the second geometry (disconnected solid skeleton).
Figure 5 shows the comparative results for the same values $\delta,
\ \mu,\ g,\ \rho_{s},\ \lambda,\, \varepsilon$ (see Table 1)  and the same
initial values.

\begin{center}
\includegraphics[width=165mm,height=60mm]{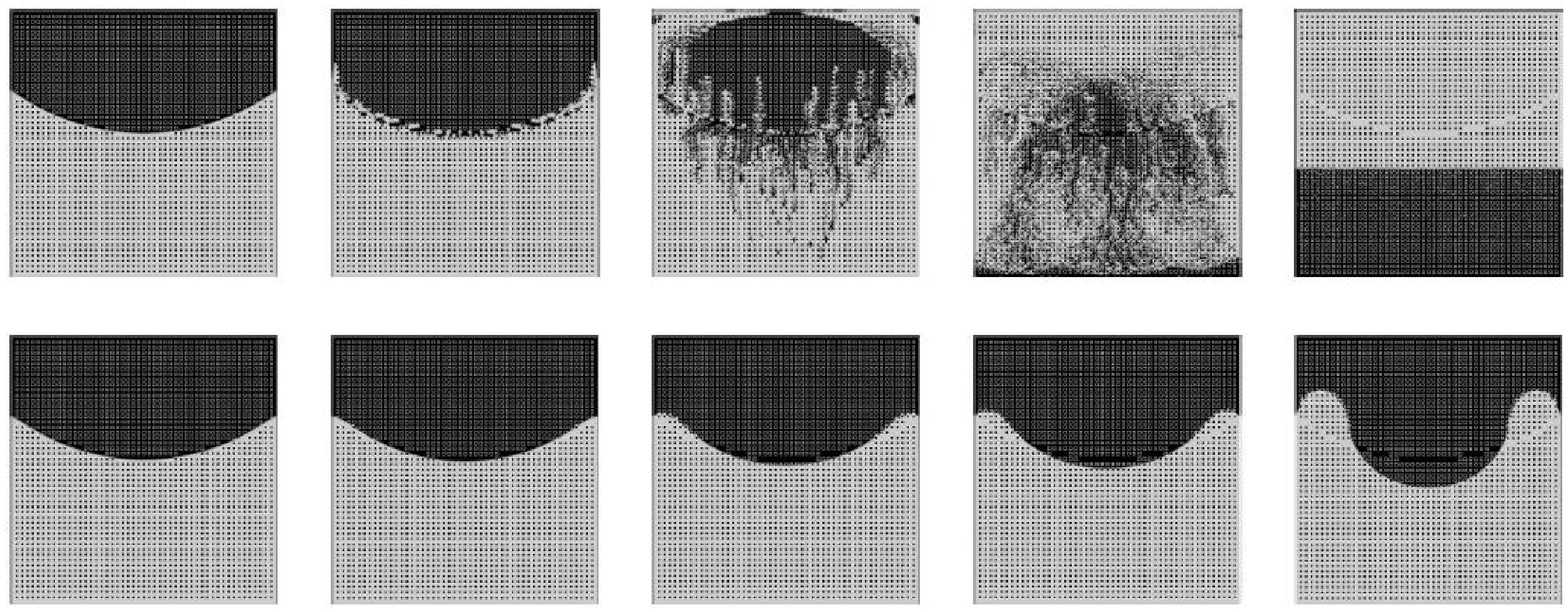}
\vskip - 0.1 cm \phantom{.} \hskip 0mm $\scriptstyle t = 50$\hskip
25mm $\scriptstyle t = 860$ \hskip 22mm $\scriptstyle t = 2631$
\hskip 22mm $\scriptstyle t = 3012$ \hskip 22mm $\scriptstyle t =
4873$
\end{center}
\begin{center}
{\scriptsize Fig.5. Comparing results of absolutely rigid and
elastic solid skeleton cases, absolutely rigid solid skeleton is
above, elastic solid skeleton is down (disconnected skeleton).}
\end{center}

\begin{center}
\begin{tabular}{|c|c|c|}
\hline
\textbf{Parameter} & \textbf{Absolutely rigid solid skeleton} & \textbf{Elastic solid skeleton} \\[5pt]
\hline
$\rho^{+}_{f}$ & $998.2$ & $998.2$\\
\hline
$\rho^{-}_{f}$ & $800$ & $800$\\
\hline
$\rho_{s}$ & $-$ & $2000$\\
\hline
$\mu^{+}$ & $10^{-2}$ & $10^{-2}$\\
\hline
$\mu^{-}$ & $9*10^{-1}$ & $9*10^{-1}$\\
\hline
$\lambda$ & $-$ & $0.5$\\
\hline
$\varepsilon$ & $2*10^{-5}$ & $2*10^{-5}$\\
\hline
$L$ & $100$ & $100$\\
\hline
\end{tabular}
\end{center}

\begin{center}
{\scriptsize Table 1. The table of values for the cases of
absolutely rigid solid skeleton and elastic solid skeleton.}
\end{center}

One may see, that for the same time of the process, 
in the Musket problem appears a mushy region (for definition see \cite{AM2}), while the  viscoelastic filtration preserves the free boundary.

As for the first geometry, we have done numerical calculations for various values $\lambda$ and $\delta$:
\begin{itemize}
 \item for $\delta = 1.25$, $\varepsilon = 2*10^{-5}$ and $\lambda \searrow 0$, there is a changing of the interface of fluids and the process of filtration is very slow;
 \item for $\delta \searrow \infty$, $\varepsilon = 2*10^{-5}$ and $\lambda = 0.5$, there is a changing of the interface of fluids and the process of filtration
become faster with increasing $\delta$;
 \item for $\delta = 1.25$, $\varepsilon = 2*10^{-5}$ and $\lambda \searrow \infty$, there is a changing of the interface of fluids and the process of filtration
become faster with increasing $\lambda$.
\end{itemize}

On the basis of numerical calculations we can conclude that

1) for the liquid motion in the absolutely rigid solid skeleton instead of the free boundary appears a mushy region, while the motion in the elastic solid skeleton preserves the free boundary;

2) when the liquids in the absolutely rigid solid skeleton completely change their positions, the liquids in the elastic solid skeleton still preserve their positions.

Therefore, the solution to the model of a viscoelastic filtration is a classical one and possesses a smooth and stable free boundary, whereas the the solution to the Musket problem will at best be only a generalized solution with mushy region instead of the free boundary.

\section{The problem statement}

\subsection{Absolutely rigid  skeleton}

We suppose that
$\Omega$ is a unit square $(0<x_{1}<1)\times(0<x_{2}<1)$, $\Omega=\Omega^{\varepsilon}_{f}\cup S^{\varepsilon}\cup\Omega^{\varepsilon}_{s}$, where
$\Omega^{\varepsilon}_{f}$ is a pore space (the domain occupied by the liquid), $\Omega^{\varepsilon}_{s}$ -- the solid skeleton and $S^{\varepsilon}$ is the  ``solid skeleton -- pore space" interface.

For the first geometry the pore space is a set of isolated capillaries
\[
\Omega^{\varepsilon}_{f}=\bigcup_{k=0}^{n-1}\big(\varepsilon k<x_{1}<\varepsilon (k+1)m\big)\times\big(0<x_{2}<1\big).
\]
Here $\varepsilon=1/n$,   $m$ is a porosity, $0<m<1$.

For the second geometry the solid skeleton is a union of disjoint squares
\[
\Omega^{\varepsilon}_{s}=\bigcup_{k=0}^{n-1}\big(\varepsilon (k+1)\beta<x_{1}<\varepsilon (k+1)(1-\beta)\big)
\times\big(\varepsilon (k+1)\beta<x_{2}<\varepsilon
(k+1)(1-\beta)\big),
\]
$2\beta=1-\sqrt{m}$.

As we have  mentioned above, two immiscible
liquids are modeled by an inhomogeneous liquid, where the density
$\rho$ can take only two constant values $\rho^{+}$ or $\rho^{-}$.
The velocity $\boldsymbol{v}$, the pressure $p_{f}$ and the density $\rho_{f}$ of the inhomogeneous liquid in the pore space $\Omega^{\varepsilon}_{f}$ are described by the Stokes system
\begin{equation}
\nabla\,\cdot\big(\mu_{1}\varepsilon^{2}\mathbb{D}(x,\boldsymbol{v})-
p_{f}\mathbb{I}\big) + \rho_{f}\,\boldsymbol{e}_{2}=0,
\,\,\,\nabla\,\cdot\boldsymbol{v}=0,
\label{2.1}
\end{equation}
and by the transport equation
\begin{equation}
\frac{\partial\rho_{f}}{\partial t}+
\boldsymbol{v}\cdot\nabla\,\rho_{f}=0,
\label{2.2}
\end{equation}
where $\boldsymbol{e}_{2}$ is the unit vector, which coincides
with the direction of the gravity.

Differential equations (2.1) and (2.2) are supplemented by the normalization condition
\begin{equation}
\int_{\Omega^{\varepsilon}_{f}}p_{f}(\boldsymbol{x},t)dx=0,
\label{2.3}
\end{equation}
the homogeneous boundary condition
\begin{equation}
\boldsymbol{v}=0, \quad\boldsymbol{x}\in S\cup S^{\varepsilon},\,\, t>0,
\label{2.4}
\end{equation}
where $S=\partial\Omega$,  and the initial condition
\begin{equation}
\rho_{f}(\boldsymbol{x},0)=\rho_{f}^{0}(\boldsymbol{x}),
\quad\boldsymbol{x}\in \Omega^{\varepsilon}_{f},
\label{2.5}
\end{equation}
\[
\rho_{f}^{0}(\boldsymbol{x})=\rho_{f}^{+}=\mbox{const}>0,\quad
\mbox{if}\,\,\boldsymbol{x}\in \Omega_{f}^{+}=
\Omega^{\varepsilon}_{f}\cap\Omega^{+},
\]
and
\[
\rho_{f}^{0}(\boldsymbol{x})=\rho_{f}^{-}=\mbox{const}>0,\quad \mbox{if}\,\,
\boldsymbol{x}\in \Omega_{f}^{-}=
\Omega^{\varepsilon}_{f}\cap\Omega^{-}.
\]

\subsection{ Elastic solid skeleton}

The joint motion of the viscous liquid in pore space and elastic solid skeleton on the microscopic level is governed  by the system (\ref{1.10}) -- (\ref{1.13}), which consists of the stationary Stokes  system
\begin{equation}
\nabla\,\cdot\big(\mu_{0}\,\mathbb{D}(x,
\frac{\partial \boldsymbol{w}_{f}}{\partial t})\big)-
\nabla p_{f} + \rho_{f}\,\boldsymbol{e}_{2}=0,
\,\,\,\nabla\,\cdot\boldsymbol{w}_{f}=0,
\label{2.6}
\end{equation}
for the displacements $\boldsymbol{w}_{f}$ and the pressure $p_{f}$ of the inhomogeneous liquid, the stationary Lam\'{e}'s equations
\begin{equation}
\nabla\cdot\big(\lambda_{0}\,\mathbb{D}(\boldsymbol{w_{s}})\big)-
\nabla p_{s}+\rho_{s}\,\boldsymbol{e}_{2}=0,\,\,\, \nabla\cdot\,\boldsymbol{w}_{s}=0,
\label{2.7}
\end{equation}
for the displacements $\boldsymbol{w}_{s}$ and the pressure $p_{s}$ of the elastic solid skeleton, and two boundary conditions on the common  ``solid skeleton -- pore space"  boundary $S^{\varepsilon}$:
\begin{equation}
\boldsymbol{w}_{f}=\boldsymbol{w}_{s},
\label{2.8}
\end{equation}
\begin{equation}
(\mu_{0}\,\mathbb{D}
(\frac{\partial \boldsymbol{w}_{f}}{\partial t})-
p_{f}\mathbb{I})\cdot\boldsymbol{n}=
(\lambda_{0}\,\mathbb{D}(\boldsymbol{w}_{s})-
p_{s}\mathbb{I})\cdot\boldsymbol{n},
\label{2.9}
\end{equation}
where $\boldsymbol{n}$ is the unit normal vector to the boundary $S^{\varepsilon}$.

The system (\ref{2.6}) -- (\ref{2.9}) is supplemented with initial and boundary conditions
\begin{equation}
\boldsymbol{w}_{f}(\boldsymbol{x},0)=0, \quad\boldsymbol{x}\in \Omega^{\varepsilon}_{f},
\label{2.10}
\end{equation}
\begin{equation}
\boldsymbol{w}=0, \quad\boldsymbol{x}\in S,\,\, t>0,
\label{2.11}
\end{equation}
where
\[
\boldsymbol{w}(x,t)=\boldsymbol{w}_{f}(x,t),\,\,\mbox{if}\,
x\in\Omega_{f},\,t>0,
\]
and
\[
\boldsymbol{w}(x,t)=\boldsymbol{w}_{s}(x,t),\,\,\mbox{if}\,
x\in\Omega_{s},\,t>0,
\]
and normalization condition
\begin{equation}
\int_{\Omega^{\varepsilon}}p(\boldsymbol{x},t)dx=0,
\label{2.12}
\end{equation}
where
\[
p(\boldsymbol{x},t)=p_{f}(\boldsymbol{x},t),\,\,\mbox{if}\,
\boldsymbol{x}\in\Omega_{f},\,t>0,
\]
and
\[
p(\boldsymbol{x},t)=p_{s}(\boldsymbol{x},t),\,\,\mbox{if}\,
\boldsymbol{x}\in\Omega_{s},\,t>0,
\]
Next, we have to complete the Cauchy problem for the transport equation  (\ref{2.2}). On the microscopic level the transport equation for the liquid density is defined only in the pore space $\Omega_{f}^{\varepsilon}$. For the exact microscopic model (see  \cite{AM}) the characteristic function $\tilde{\chi}^{\varepsilon}$ of the pore space is an unknown function and defined as a solution to the Cauchy problem
\[
\frac{\partial\tilde{\chi}^{\varepsilon}}{\partial t}+\frac{\partial\boldsymbol w}{\partial t}\cdot\,
\nabla\tilde{\chi}^{\varepsilon}=0, \quad
\tilde{\chi}^{\varepsilon}(\boldsymbol{x},0)=
\chi^{\varepsilon}(\boldsymbol{x}).
\]
It means that the solid -- liquid interface $S^{\varepsilon}$ in the exact model is a material surface and we do not need boundary conditions for the liquid density on  $S^{\varepsilon}$. But our basic dynamic system here is a linear one, where $\tilde{\chi}^{\varepsilon}=\chi^{\varepsilon}$ is the given function. Therefore the solid -- liquid interface $S^{\varepsilon}$ is no longer the material surface and we need the boundary condition for the liquid density on the part of $S^{\varepsilon}$, where the liquid "enters" into the pore space. To avoid this, we extend the Cauchy problem (\ref{2.2}) onto hole domain $\Omega$.  At first, we suppose that the function $\rho_{f}^{0}(\boldsymbol{x})$ is defined in hole domain $\Omega$ and
\[
\rho_{f}^{0}(\boldsymbol{x})=\rho_{f}^{+} \,\,\mbox{in}\,\,\Omega^{+},\,\,
\rho_{f}^{0}(\boldsymbol{x})=\rho_{f}^{-} \,\,\mbox{in}\,\,\Omega^{-}.
\]
Finally, we rewrite  the Cauchy problem for the liquid density as the Cauchy problem for the density of a mixture
\[
\rho=\chi^{\varepsilon}\rho_{f}+(1-\chi^{\varepsilon})\rho_{s}
\]
in the form
\begin{equation}
\frac{\partial\rho}{\partial t}+\frac{\partial \boldsymbol{w}}{\partial t}\cdot\,\nabla\rho=0,\,\,
(\boldsymbol{x},t)\in\Omega_{T}; \quad
\rho(\boldsymbol{x},0)=\rho^{0}(\boldsymbol{x}),
\,\,\boldsymbol{x}\in\Omega,
\label{2.13}
\end{equation}
where
\[
\rho^{(0)}(\boldsymbol{x})=
\chi^{\varepsilon}\rho_{f}^{0}(\boldsymbol{x})+
(1-\chi^{\varepsilon})\rho_{s}.
\]
\section{The computational algorithm}

\subsection{The absolutely rigid solid skeleton}

As the numerical method for computer simulation (\ref{2.1}) --
(\ref{2.2}) system of equations, supplemented with appropriate
conditions (\ref{2.3}) -- (\ref{2.5}), we chose the large-particle
method. This method is based on splitting the original differential
equations in accordance with the physical processes they represent.
The large-particle method is the development of Harlow's method of
"particle-in-cell", which refers to the technique used to solve a
certain class of partial differential equations, including
the Rayleigh-Taylor instability for the Musket free boundary problem.

The solution process for an evolution system (\ref{2.1}) --
(\ref{2.2}) is very difficult, because we have the presence of
stratification. Thus the heterogeneity of the fluid requires
additional calculation of the density field. To find the unknown
functions the calculation process can be represented as a three-step
scheme.\\

\textbf{Step 1}: One has to calculate intermediate velocity $\boldsymbol{\tilde{v}}$ from the equation

\begin{equation}
\nabla\,\cdot\big(\mu_{1}\varepsilon^{2}
\mathbb{D}(\boldsymbol{\tilde{v}})\big) + \rho
\mathbf{e}_{2} = 0,\label{3.1}
\end{equation}
and to find intermediate density $\tilde{\rho}$ from the equation

\begin{equation}
\frac{\partial \tilde{\rho}}{\partial t} + \boldsymbol{\tilde{v}}\,\cdot \nabla\rho = 0,\label{3.2}
\end{equation}

\textbf{Step 2}: The main difficulty in the numerical solution of
equations (2.1), (2.2) associated with the calculation of the
pressure field. The first significant success in overcoming these
difficulties has been achieved through the idea of artificial
compressibility \cite{Galtsev:Article7}. It is very important to
calculate unknown $p$ from the equation:

\begin{equation}
\frac{\partial p}{\partial t}+c_{p,f}^{2} \nabla \,\cdot \boldsymbol{\tilde{v}}=0.\label{3.3}
\end{equation}

Here we need to satisfy the solenoidal condition ($\nabla\,\cdot
\boldsymbol{\tilde{v}} = 0$). It is necessary to increase $c_{p}$
each time step, while it will be equal to given accuracy.\\

\textbf{Step 3}: We calculate the final values of the velocity $\boldsymbol{v}$ from the
Stokes equation on the next time step:

\begin{equation}
\nabla\,\cdot\big(\mu_{1}\varepsilon^{2}\mathbb{D}(\boldsymbol{v})\big) -
\nabla p + \tilde{\rho} \mathbf{e}_{2} = 0 \label{3.4}
\end{equation}
and the density $\rho$ from the equation:

\begin{equation}
\frac{\partial {\rho}}{\partial t} + \boldsymbol{v}\,\cdot
\nabla\tilde{\rho} = 0,\label{3.5}
\end{equation}

So, let us  write the basic finite-difference scheme.
The domain of integration is covered by
a stationary (Eulerian) difference grid of arbitrary form (to
abbreviate the exposition, a rectangular grid in a
planar domain is considered, see Fig.6):
\[
\Omega_{f} =\begin{pmatrix}
x_{1}^{(i+1/2)} = ih_{1}, \,\, h_{1}>0;\,\, i=0,1,...,N_{1};\\
x_{2}^{(j+1/2)} = jh_{2}, \,\, h_{2}>0;\,\, j=0,1,...,N_{2};
\end{pmatrix}
\]
where $h_{1}$, $h_{2}$ is the size of the grid, $N_{1}$, $N_{2}$ are
the numbers of grid cells, respectively, in the $x_{1}$ and $x_{2}$
direction (the point with coordinates $(i,j)$ matches with the
center of the cell). Here, as in the original splitting method, we
use the "checkerboard" grid. This makes it possible  clearly
interpret each cell as element of volume, which is characterized by
a calculated pressure and density in it's center.

\begin{center}
\includegraphics[width=95mm,height=75mm]{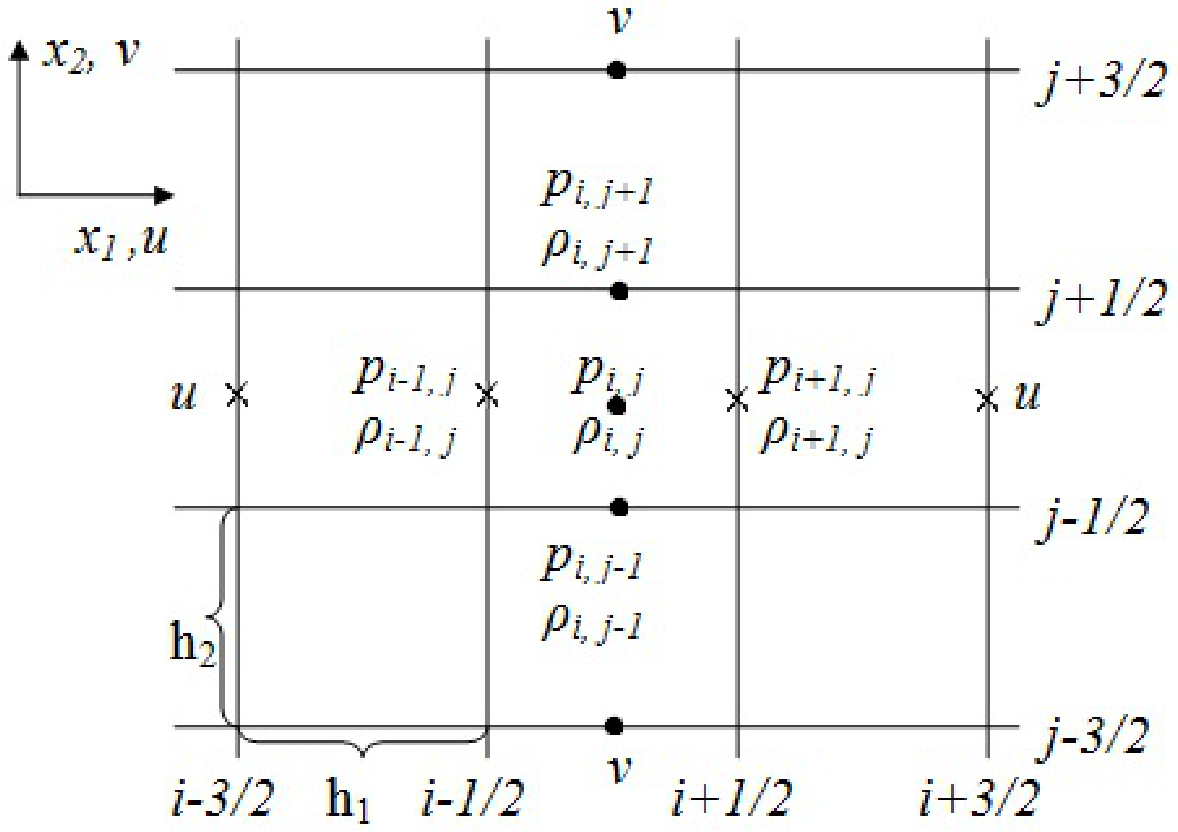}
\end{center}
\begin{center}
\phantom{.}\vskip -0.7 cm $\phantom{.}$ {\scriptsize  Fig.6.
Template of the grid.}
\end{center}

Only quantities relating to the cell as a whole are varied, while the fluid is assumed to be momentarily
constrained. For this reason, convective terms, corresponding to
translational effects, are dropped from the equations. In the
remaining equations, $\rho$ is brought forward outside the
differentiation symbol, and the equations (2.1)--(2.2) are solved for
the derivatives of $u$, $v$ with respect to the time.

An elementary finite-difference approximation of the equation (3.1)
yields the following expressions:

\[
\frac{\tilde{u}_{i+3/2,j}^{n}-2\,\tilde{u}_{i+1/2,j}^{n}+\tilde{u}_{i-/2,j}^{n}}{h_{1}^{2}}
= - \frac{\tilde{u}_{i+1/2,j+1}^{n}-2\,\tilde{u}_{i+1/2,j}^{n}+\tilde{u}_{i+1/2,j-1}^{n}}{h_{2}^2},
\]
\[
\frac{\tilde{v}_{i+1,j+1/2}^{n}-2\,\tilde{v}_{i,j+1/2}^{n}+\tilde{v}_{i-1,j+1/2}^{n}}{h_{1}^{2}}
= - \frac{\tilde{v}_{i,j+3/2}^{n}-2\,\tilde{v}_{i,j+1/2}^{n}+\tilde{v}_{i,j-1/2}^{n}}{h_{2}^2}
- \frac{\rho_{i,j+1/2}}{\mu_{1}\varepsilon^{2}},
\]
where $n$ -- number of the time step, $\tilde{u}$, $\tilde{v}$ -- are intermediate values of the flow velocity. \\
Quantities with fractional subscripts relate to cell boundaries,
e.g.:

$$\tilde{v}_{i+1/2, j} = \frac{\tilde{v}_{i,j}+\tilde{v}_{i+1,j}}{2},$$

As in the large particles method, we need to calculate the flux density (3.2)
through the cell boundaries. It is
assumed throughout that the density of the large particle is in the
motion only owing to the velocity component normal to the boundary.
This values of the flow parameters in the next
time layer are computed according to the following formulas (the
direction of the flow is from left to right and upwards):

\begin{equation}
\mu_{1}\varepsilon^{2} \frac{\tilde{\rho}_{i,j}^{n+1}-\tilde{\rho}_{i,j}^{n}}{\tau} =
-\frac{(\rho \tilde{u})^{n+1}_{i+1/2,j}-(\rho \tilde{u})^{n+1}_{i-1/2,j}}{h_{1}} -
\frac{(\rho \tilde{v})^{n+1}_{i,j+1/2}-(\rho
\tilde{v})^{n+1}_{i,j-1/2}}{h_{2}}.\label{3.4}
\end{equation}

In cases, where we need to determine the function values in the grid
points, not to meet their terms in Fig. 6, we used the
arithmetic average, for example: $$\rho_{i+1/2} =
\frac{1}{2}(\rho_{i+1,j} + \rho_{i,j}).$$

At the final step (when a discrete model of the medium is being
used) one should carry out an additional calculation of the density;
this smoothness out fluctuations and increases the accuracy of the
computations. Combining  different representations of the steps,
one obtains a series of difference schemes; this series, which
constitutes the large-particle method.

Chosen method may be interpreted from various points of
view: the splitting method, the mixed Euler–Lagrange method,
computation in local Lagrangian coordinates with
scaling on the previous grid, the difference notation for
conservation laws for a fluid element (large
particle), and the Eulerian difference scheme.

When we replace the differential problem by the finite-difference
representation, we should pay special attention to the approximation
of boundary conditions, because their specific approximation affects the
correctness of the method and stability of the scheme, as well as the
velocity of convergence.

The boundary conditions are formulated by introducing series of
fictitious cells (so that every computation point becomes an
interior point and the same algorithm is maintained for all cells).
One layer is sufficient for the first-order approximation scheme,
two layers are sufficient for the second order. As the result, in
the case, where side walls are the solid surface, the impermeability
condition is represented as
\begin{equation}
u_{-1/2, j} = 0,\label{3.5}
\end{equation}
and the slip condition as
\begin{equation}
v_{i-1/2, j+1/2} = 0.\label{3.6}
\end{equation}

In the planar case, the geometrical characteristics of fractional
cells may be determined by direct measurement. In the
axially-symmetric case one needs an additional computation,
incorporating the distance of each fractional cell from the axis of
symmetry. The difference formulas for fractional cells are obtained
by a slight modification of the difference formulas for whole cells.

\subsection{ The elastic solid skeleton}

The computer simulation of joint motion of the viscous liquid (in
the pore space) and elastic solid skeleton (at the microscopic
level) is the numerical solution of the system, which consists of
the stationary Stokes system (\ref{2.6}) for the displacement
$\boldsymbol{w}_{f}$ and the pressure $p_{f}$ of the inhomogeneous
liquid and the stationary Lam\'{e}'s equations (\ref{2.7}) for the
displacement $\boldsymbol{w}_{s}$ and the pressure $p_{s}$ of the
elastic solid skeleton with appropriate conditions on the common
border.

As the numerical method to simulate the Rayleigh-Taylor instability
with an elastic component of the skeleton, we chose the same method,
which has been described above (large particles method).Considerate
domain is replaced by the system of fluid particles, which match
with the cell of Eulerian mesh (Fig.6).

The suggested algorithm consists of four main steps:

\textbf{Step 1}. One has to solve the system of Lam\'{e}'s equations
(\ref{2.7}) with given boundary and initial conditions for the
displacement $\boldsymbol {w}_{s}$ and the solid pressure  $p_{s}$:
\begin{equation}
\boldsymbol{w}_{s}=\boldsymbol{w}_{f},\,\,\boldsymbol{x}\in S^{\varepsilon},\,\,
\boldsymbol{w}_{s}=0,\,\,\boldsymbol{x}\in S.
\label{3.7}
\end{equation}

The solution of this system is no different from the solution of the
Stokes system (Subsection 3.1). We enter the artificial
compressibility and calculate unknown $p_{s}$ from the equation::

\begin{equation}
\frac{\partial p_{s}}{\partial t}+c_{p,s}^{2} \nabla \,\cdot
\boldsymbol{w_{s}}=0.\label{3.8}
\end{equation}

It is necessary to increase $c_{p,s}$ each time step, while the
solenoidal condition ($\nabla\,\cdot \boldsymbol{w}_{s} = 0$) will
be satisfied, since the fluid is incompressible.

Also it should be noted, that the procedure for two-dimensional
difference scheme of the Lam\'{e}'s equation is identical to the
previous. And scheme itself is easier, unlike for Stokes equation,
because we needn't to introduce the necessary differential
approximation for the velocity.

An elementary finite-difference approximation of the equation (2.7)
yields the following expressions:

\[
\frac{w_{1\,i+3/2,j}^{n}-2\,w_{1\,i+1/2,j}^{n}+w_{1\,i-/2,j}^{n}}{h_{1}^{2}}
= - \frac{w_{1\,i+1/2,j+1}^{n}-2\,w_{1\,i+1/2,j}^{n}+w_{1\,i+1/2,j-1}^{n}}{h_{2}^2},
\]
\[
\frac{w_{2\,i+1,j+1/2}^{n}-2\,w_{2\,i,j+1/2}^{n}+w_{2\,i-1,j+1/2}^{n}}{h_{1}^{2}}
= - \frac{w_{2\,i,j+3/2}^{n}-2\,w_{2\,i,j+1/2}^{n}+w_{2\,i,j-1/2}^{n}}{h_{2}^2}
- \frac{\rho_{s\,i,j+1/2}}{\lambda_{0}},
\]
where $n$ -- number of the time step, $w_{1}$, $w_{2}$ are
displacement values, $\rho_{s}$ is the density of the elastic solid skeleton. \\

\textbf{Step 2}. Using known values $\boldsymbol {w}_{s}$ and
$p_{s}$, find normal stress on the boundary $S^{\varepsilon}$:

\begin{equation}
(\lambda_{0} \mathbb{D}(\boldsymbol{
w_{s}})-p_{s}\mathbb{I})\boldsymbol{n}=\boldsymbol{A},\label{3.9}
\end{equation}
where $\boldsymbol {n}$ is the unit normal to the boundary
$S^{\varepsilon}$.\\

\textbf{Step 3}. Then solve the system of Stokes equations
(\ref{2.6}) in $\Omega_ {f}^{\varepsilon}$ just as for absolutely
rigid solid skeleton, repeating Step 1 -- Step 3 from the previous
subsection, with the condition on the common boundary
$S^{\varepsilon}$:

\begin{equation}
(\mu_{0}\mathbb{D}(\frac{\partial\boldsymbol w_{f}}{\partial
t})-p_{f}\mathbb{I})\boldsymbol{n} = (\lambda_{0}
\mathbb{D}(\boldsymbol{w_{s}})-p_{s}\mathbb{I})\boldsymbol{n},\label{3.10}
\end{equation}
where we know the right side of the equality from the previous step.

This boundary condition suggests, that the vector of displacement
and pressure satisfies continuity of normal stresses on the common
boundary between liquid and elastic solid skeleton.

The found value of the fluid velocity replace to the transport
equation (\ref{2.13}) in $\Omega_{f}^{\varepsilon}$, where we find
the density value $\rho_{f}$ for the next time step.\\

\textbf{Stage 4}. When the velocity $\partial\boldsymbol
w_{f}/\partial t$ is known, in the liquid part we determine
$\boldsymbol w_{s}$ on the next time step from the continuity of
normal stresses and condition (2.8).

Thus, considering behavior of fluids on the boundary of solid
skeleton and, solving the system of equations (\ref{2.6}),
(\ref{2.7}), (\ref{2.13}) with appropriate initial and boundary
conditions (\ref{2.8}) -- (\ref{2.12}), we obtain the numerical
approximation of joint motion of fluid and elastic solid skeleton.

\begin{center} \textbf{Conclusions}
\end{center}

In the present  paper we have shown, how  to model physical processes using modern methods of the mathematical analysis.  We started with the free boundary problem for a joint motion of two immiscible incompressible fluids on the   microscopic level. Theoretically  this mathematical model is the most suitable  model, describing the given physical process. But this model has no practical  value, because we have to solve the problem in the physical domain of some hundreds meters, while the coefficients oscillate on the physical size in some microns. The practical  value of the model appears only after homogenization.  In turn, the homogenization has at least two levels of approximation, which depend on the dimensionless criteria of the physical problem. The first level of approximation is the well -- known Muskat problem.  The second level of approximation of the free boundary problem on the  microscopic level is the Muskat problem for a viscoelastic filtration. In our numerical calculations for the periodic structure, we simulated the homogenization  by increasing the number of elementary cells per unit volume. The numerical results show that the solution to the Muskat problem is unstable, the free boundary louses its sharp structure and transforms into some domain (mushy region), occupied by the mixture of two liquids. Whereas, the solution to the Muskat problem for a viscoelastic filtration remains a classical one with a smooth free boundary.  That is, the Muskat problem for a viscoelastic filtration is a natural generalization of the classical Muskat problem, which still remains unsolved as a mathematical task, and extremely  difficult for the numerical realization.

\begin{center} \textbf{Acknowledgment}
\end{center}

This research  is partially supported  by the Federal  Program "Research and scientific-pedagogical brainpower of  Innovative Russia" for 2009-2013 (State Contract  02.740.11.0613).

\end{document}